\documentclass[reqno]{amsart}
\usepackage{}
\usepackage{stmaryrd}
\usepackage{mathrsfs}
\usepackage{cases}
\usepackage{amsfonts}
\usepackage{amssymb}
\usepackage{amsmath}
\usepackage{tikz}
\newtheorem{theorem}{Theorem}[section]

\newtheorem{proposition}[theorem]{Proposition}

\theoremstyle{definition}
\newtheorem{definition}[theorem]{Definition}

\newtheorem{remark}[theorem]{Remark}
\numberwithin{equation}{section}

\newcommand{\blankbox}[2]

\let\al=\alpha
\let\a=\alpha
\let\b=\beta

\let\r=\rho

\let\f=\frac

\let\Om=\Omega

\let\ep=\epsilon

\def\bbN{\mathbb{N}}

\newcommand{\be}{\begin{equation*}}
\newcommand{\ee}{\end{equation*}}
\newcommand{\ben}{\begin{equation}}
\newcommand{\een}{\end{equation}}
\newcommand{\bn}{\begin{enumerate}}
\newcommand{\en}{\end{enumerate}}
\newcommand{\bs}{\backslash}
\newcommand{\fy}{\infty}

\def\rr{{\mathbb R}}
\def\rn{{{\rr}^n}}

\begin{document}
\title[Limiting weak-type behaviors with rough kernel]
{Limiting weak-type behaviors for factional maximal operators and fractional integrals
with rough kernel}

\author{GUOPING ZHAO}
\address{School of Applied Mathematics, Xiamen University of Technology, Xiamen, 361024, P.R.China}
\email{guopingzhaomath@gmail.com}
\author{WEICHAO GUO}
\address{School of Science, Jimei University, Xiamen, 361021, P.R.China}
\email{weichaoguomath@gmail.com}

\begin{abstract}
By a reduction method, the limiting weak-type behaviors of factional maximal operators and fractional integrals
are established without any smoothness assumption on the kernel, which essentially improve and extend previous results.
As a byproduct, we characterize the boundedness of several operators by the membership of their kernel in Lebesgue space on sphere.
\end{abstract}



\thanks{Supported by the NSF of China (Nos.11701112, 11771388, 11671414, 11601456).}

\maketitle

\section{Introduction}
In order to study the lower bounds for the weak-type constants of singular integral operators $T_{\Omega}$ with homogeneous kernel $\Om$, with some regularity conditions on $\Om$,
the limiting weak-type behaviors were first established in \cite{Jan_Trans_2006}.
Then, with weaker regularity conditions on $\Om$,
Ding-Lai establish the same limiting weak-type behaviors of $T_{\Omega}$ in \cite{DL2},
where they also obtain the limiting weak-type behaviors of fractional integral operators $T_{\Omega}^{\alpha}$
with some regularity conditions on $\Om$.
At the same time,
the corresponding results of maximal and fractional maximal operators with homogeneous kernels were
established in \cite{DL1}.
Recently, a refined limiting weak-type behavior associated with stronger convergence, namely, type-1 and type-2 convergence (see Proposition \ref{proposition, comparison for limiting behavior}), was established by Guo-He-Wu \cite{GHW}.
For a study in vector-valued setting, we refer to \cite{HGW-JMAA}.
One can see \cite{Tang-Wang_NA2020} for the limiting weak-type behavior for
multilinear fractional integrals.

Let $\alpha\in (0,n)$ and $\Omega$ be a homogeneous function of degree zero.
The integral operator $T_{\Omega}^{\alpha}$ with homogeneous kernel $\Omega$ is defined as
\be
  T_{\Omega}^{\alpha}f(x):=\int_{\mathbb{R}^n}\frac{\Omega(x-y)}{|x-y|^{n-\alpha}}f(y)dy,
\ee
and the corresponding fractional maximal operator is defined by
\be
  M_{\Omega}^{\alpha}f(x):=\sup_{r>0}\frac{1}{r^{n-\alpha}}\int_{B(x,r)}|\Omega(x-y)f(y)|dy.
\ee

In the previous works, in order to establish the limiting weak-type behavior
for factional integral, some angular regularity must be added to the kernel $\Om$.
More exactly, the definition of angular regularity, namely, the Dini-condition is as following.

\begin{definition}[$L^q_s$-Dini condition]\label{definition, Dini condition}
  Suppose $\Omega$ is a homogeneous function of degree zero. Let $1\leq q\leq \infty, 0\leq s<n$. We say that
  $\Omega$ satisfies $L^q_{s}$-Dini condition if
  \begin{enumerate}
    \item $\Omega\in L^q(\mathbb{S}^{n-1})$,
    \item $\int_0^1 \frac{{\omega}_q(t)}{t^{1+s}}dt<\infty$,
  \end{enumerate}
  where ${\omega}_q$ is called the (modified) integral continuous modulus of $\Omega$ of degree $q$, defined by
  \begin{equation}
    {\omega}_q(t):=\left(\sup_{|h|\leq t}\int_{\mathbb{S}^{n-1}}|\Omega(x'+h)-\Omega(x')|^qd\sigma(x')\right)^{1/q}.
  \end{equation}
\end{definition}
In order to compare different types of limiting weak-type behaviors, we recall a proposition established in \cite{GHW}.
\begin{proposition}[see \cite{GHW}]\label{proposition, comparison for limiting behavior}
  Let $0<p<\infty$. Suppose that $f\in L^{p,\infty}(\mathbb{R}^n)$, and
  $|\{x\in \mathbb{R}^n: |f(x)|=\lambda\}|=0$ for all $\lambda > 0$.
  Let $\{f_{(t)}\}_{t>0}$ be a sequence of measurable functions.
  Then for the following three statements:
  \begin{enumerate}
  \item $\forall\varepsilon>0$, $\exists A_{\varepsilon}\subset \mathbb{R}^n$, s.t.,
  $|A_{\varepsilon}|<\varepsilon$ and $\displaystyle\lim_{t\rightarrow 0^{+}}\|f-f_{(t)}\|_{L^{p,\infty}(\mathbb{R}^n\backslash A_{\varepsilon})}=0$;
  \item $\displaystyle\lim_{t\rightarrow 0^{+}}|\{x\in \mathbb{R}^n:\, |f_{(t)}(x)-f(x)|>\lambda\}|=0$,\, $\forall\,\lambda>0$;
  \item $\displaystyle\lim_{t\rightarrow 0^{+}}|\{x\in \mathbb{R}^n: \,|f_{(t)}(x)|>\lambda\}|=|\{x\in \mathbb{R}^n: |f(x)|>\lambda\}|$,\, $\forall\,\lambda>0$;
\end{enumerate}
we have
\begin{equation*}
  (1)\Rightarrow (2) \Rightarrow (3),\ (3)\nRightarrow (2),\ (2)\nRightarrow (1).
\end{equation*}
\end{proposition}
We say that a sequence of functions $f_{(t)}$ tends to $f$
  in the sense of type-$i$,
  if the statement $(i)$ in Proposition \ref{proposition, comparison for limiting behavior} is valid, $i=1,2,3$.

Based on our limited knowledge for the limiting weak-type behaviors of fractional maximal and fractional integrals,
the best result so far comes from \cite{GHW} as we shall list as following.
Denote $f_t(\cdot):=\frac{1}{t^n}f(\f{\cdot}{t})$.

\medskip

\textbf{Theorem A} (cf. \cite[Theorem 1.9,1.10]{GHW}).\quad
Let $\alpha\in (0,n)$, $f\in L^1(\rn)$ be a nonnegative function.
Suppose that $\Omega$ is a homogeneous function of degree zero
and the maximal operator $M_{\Omega}^{\alpha}$ is bounded from $L^1$ to $L^{\frac{n}{n-\alpha},\infty}$.
If $\Om$ satisfies the $L^1_{\alpha}$-Dini condition, we have
\begin{enumerate}
\item    $\Omega\in L^{\frac{n}{n-\alpha}}(\mathbb{S}^{n-1})$, $\|\Omega \|_{L^{\frac{n}{n-\alpha}}(\mathbb{S}^{n-1})}
\sim
\Big\|\frac{\Omega(\cdot)}{|\cdot|^{n-\alpha}}\Big\|_{L^{\frac{n}{n-\alpha},\infty}}
\lesssim \|M_{\Omega}^{\alpha}\|_{L^1\rightarrow L^{\frac{n}{n-\alpha},\infty}}$;
\item   $\displaystyle\lim_{t\rightarrow 0^{+}}\Big|\Big\{x\in \mathbb{R}^n: \Big|M^{\alpha}_{\Omega}f_t(\cdot)-\frac{|\Omega(\cdot)|}{|\cdot|^{\frac{n}{n-\alpha}}}\|f\|_{L^1(\rn)}\Big|>\lambda\Big\}\Big|=0$, \, $\forall\,\lambda>0$.
\end{enumerate}
If $\Om$ satisfies the $L^{\frac{n}{n-\alpha}}_0$-Dini condition, we have
\begin{enumerate}
\item    $\Omega\in L^{\frac{n}{n-\alpha}}(\mathbb{S}^{n-1})$, $\|\Omega \|_{L^{\frac{n}{n-\alpha}}(\mathbb{S}^{n-1})}
\sim
\Big\|\frac{\Omega(\cdot)}{|\cdot|^{n-\alpha}}\Big\|_{L^{\frac{n}{n-\alpha},\infty}}
\lesssim \|M_{\Omega}^{\alpha}\|_{L^1\rightarrow L^{\frac{n}{n-\alpha},\infty}}$;
\item    $\displaystyle\lim_{t\rightarrow 0^{+}} \Big\|M^{\alpha}_{\Omega}f_t(\cdot)-\frac{|\Omega(\cdot)|}{|\cdot|^{n-\alpha}}\|f\|_{L^1(\rn)}\Big\|_{L^{\frac{n}{n-\alpha},\infty}(\mathbb{R}^n\backslash B(0,\rho))}=0$
    for every $\rho>0$.
\end{enumerate}

In short, under Dini-condition assumptions on the kernel $\Om$,
Theorem A establishes the limiting weak-type behaviors for $M_{\Om}^{\al}$ in the sense of type-1 and type-2, respectively.

Next, we list the corresponding results for the fractional integral operator in \cite{GHW}.

\medskip

\textbf{Theorem B} (cf. \cite[Theorem 1.12,1.13]{GHW}).\quad
Let $\alpha\in (0,n)$, $f\in L^1(\rn)$ be a nonnegative function.
Suppose that $\Omega$ is a homogeneous function of degree zero and
the factional integral operator $T_{\Omega}^{\alpha}$ is bounded from $L^1$ to $L^{\frac{n}{n-\alpha},\infty}$.
If $\Om$ satisfies the $L^1_{\alpha}$-Dini condition, we have
\begin{enumerate}
\item    $\Omega\in L^{\frac{n}{n-\alpha}}(\mathbb{S}^{n-1})$, $\|\Omega \|_{L^{\frac{n}{n-\alpha}}(\mathbb{S}^{n-1})}
\sim
\Big\|\frac{\Omega(\cdot)}{|\cdot|^{n-\alpha}}\Big\|_{L^{\frac{n}{n-\alpha},\infty}}
\lesssim \|T_{\Omega}^{\alpha}\|_{L^1\rightarrow L^{\frac{n}{n-\alpha},\infty}}$;
\item   $\displaystyle\lim_{t\rightarrow 0^{+}}\Big|\Big\{x\in \mathbb{R}^n: \Big|T^{\alpha}_{\Omega}f_t(\cdot)-\frac{\Omega(\cdot)}{|\cdot|^{\frac{n}{n-\alpha}}}\|f\|_{L^1(\rn)}\Big|>\lambda\Big\}\Big|=0$, $\forall\,\lambda>0$.
\end{enumerate}
Moreover, if $\Om$ satisfies the $L^{\frac{n}{n-\alpha}}_0$-Dini condition, we have
\begin{enumerate}
\item    $\Omega\in L^{\frac{n}{n-\alpha}}(\mathbb{S}^{n-1})$, $\|\Omega \|_{L^{\frac{n}{n-\alpha}}(\mathbb{S}^{n-1})}
\sim
\Big\|\frac{\Omega(\cdot)}{|\cdot|^{n-\alpha}}\Big\|_{L^{\frac{n}{n-\alpha},\infty}}
\lesssim \|T_{\Omega}^{\alpha}\|_{L^1\rightarrow L^{\frac{n}{n-\alpha},\infty}}$;
\item    $\displaystyle\lim_{t\rightarrow 0^{+}} \Big\|T^{\alpha}_{\Omega}f_t(\cdot)-\frac{\Omega(\cdot)}{|\cdot|^{n-\alpha}}\|f\|_{L^1(\rn)}\Big\|_{L^{\frac{n}{n-\alpha},\infty}(\mathbb{R}^n\backslash B(0,\rho))}=0$
    for every $\rho>0$.
\end{enumerate}

Under the Dini-condition of $\Om$,
Theorem B establishes the limiting weak-type behaviors for $T_{\Om}^{\al}$ in the sense of type-1 and type-2, respectively.

In this paper, our first main goal is to establish the limiting weak-type behaviors in the sense of  type-1 for $M_{\Om}^{\al}$ and $T_{|\Om|}^{\al}$
without any smoothness assumption on $\Om$, which essentially improves Theorem A.
Surprisingly, limiting weak-type behaviors can help us in turn to
establish several equivalent characterizations associated with the boundedness of $M_{\Omega}^{\alpha}$ and $T_{|\Omega|}^{\alpha}$.
The following is our first main result:
\begin{theorem}\label{thm, maximal}
Let $\alpha\in (0,n)$, $1<r<\infty$.
Suppose that $\Omega\in L^1(\mathbb{S}^{n-1})$ is a homogeneous function of degree zero.
The following statements are equivalent.
\begin{enumerate}
  \item    $\Omega\in L^{\frac{n}{n-\alpha}}(\mathbb{S}^{n-1})$;
  \item   $M_{\Omega}^{\alpha}$ is bounded from $L^1$ to $L^{\frac{n}{n-\alpha},\infty}$;
  \item   $T_{|\Omega|}^{\alpha}$ is bounded from $L^1$ to $L^{\frac{n}{n-\alpha},\infty}$;
  \item   $M_{\Omega}^{\alpha}$ is bounded from $L^1(l^r)$ to $L^{\frac{n}{n-\alpha},\infty}(l^r)$;
  \item   $T_{|\Omega|}^{\alpha}$ is bounded from $L^1(l^r)$ to $L^{\frac{n}{n-\alpha},\infty}(l^r)$.
\end{enumerate}
Moreover, let $f\in L^1(\rn)$ and $\{f_j\}_{j\in\bbN}\in L^1(l^r)$, if one of the above statements holds, we have
\begin{enumerate}
  \item[(a)]
   $\|\Omega \|_{L^{\frac{n}{n-\alpha}}(\mathbb{S}^{n-1})}
   \sim
   \Big\|\frac{\Omega(\cdot)}{|\cdot|^{n-\alpha}}\Big\|_{L^{\frac{n}{n-\alpha},\infty}}
   \sim \|M_{\Omega}^{\alpha}\|_{L^1\rightarrow L^{\frac{n}{n-\alpha},\infty}}
   \sim \|T_{|\Omega|}^{\alpha}\|_{L^1\rightarrow L^{\frac{n}{n-\alpha},\infty}} \\
   \sim \|M^{\alpha}_{\Om}\|_{L^1(l^r)\rightarrow L^{\frac{n}{n-\alpha},\infty}(l^r)}
   \sim \|T^{\alpha}_{|\Om|}\|_{L^1(l^r)\rightarrow L^{\frac{n}{n-\alpha},\infty}(l^r)}$;
  \item[(b)]
    $\displaystyle\lim_{t\rightarrow 0^{+}} \Big\|M^{\alpha}_{\Omega}f_t(\cdot)-\frac{|\Omega(\cdot)|}{|\cdot|^{n-\alpha}}\|f\|_{L^1(\rn)}\Big\|_{L^{\frac{n}{n-\alpha},\infty}(\mathbb{R}^n\backslash B(0,\rho))}=0$
    for every $\rho>0$;
  \item[(c)]
    $\displaystyle\lim_{t\rightarrow 0^{+}} \Big\|T^{\alpha}_{|\Omega|}f_t(\cdot)-\frac{|\Omega(\cdot)|}{|\cdot|^{n-\alpha}}\|f\|_{L^1(\rn)}\Big\|_{L^{\frac{n}{n-\alpha},\infty}(\mathbb{R}^n\backslash B(0,\rho))}=0$
    for every $\rho>0$;
  \item[(d)]
    $\displaystyle\lim_{t\rightarrow 0^{+}}
    \bigg\|
     \Big\|
      \big\{
       M^{\alpha}_{\Omega}f_{j,t}(\cdot)-\frac{|\Omega(\cdot)|}{|\cdot|^{n-\alpha}}\|f_j\|_{L^1(\rn)}\big
      \}_{j\in\bbN}
     \Big\|_{l^r}
    \bigg\|_{L^{\frac{n}{n-\alpha},\infty}(\mathbb{R}^n\backslash B(0,\rho))}
    =0$
    for every $\rho>0$;
  \item[(e)]
    $\displaystyle\lim_{t\rightarrow 0^{+}}
    \bigg\|
     \Big\|
      \big\{
       T^{\alpha}_{|\Omega|}f_{j,t}(\cdot)-\frac{|\Omega(\cdot)|}{|\cdot|^{n-\alpha}}\|f_j\|_{L^1(\rn)}\big
      \}_{j\in\bbN}
     \Big\|_{l^r}
    \bigg\|_{L^{\frac{n}{n-\alpha},\infty}(\mathbb{R}^n\backslash B(0,\rho))}
    =0$
    for every $\rho>0$.
\end{enumerate}
\end{theorem}

For the fractional operator without any smoothness assumption, we give the limiting weak-type behaviors in the sense of  type-1  as following,
which essentially improve Theorem B.

\begin{theorem}\label{thm, fractional}
Let $\alpha\in (0,n)$, $f\in L^1(\rn)$.
Suppose that $\Omega\in L^{\frac{n}{n-\al}}(\mathbb{S}^{n-1})$ is a homogeneous function of degree zero.
Then $T_{\Omega}^{\alpha}$ is bounded from $L^1$ to $L^{\frac{n}{n-\alpha},\infty}$ and from $L^1(l^r)$ to $L^{\frac{n}{n-\alpha},\infty}(l^r)$, and
\begin{enumerate}
  \item
  $\|\Omega \|_{L^{\frac{n}{n-\alpha}}(\mathbb{S}^{n-1})}
   \sim
   \Big\|\frac{\Omega(\cdot)}{|\cdot|^{n-\alpha}}\Big\|_{L^{\frac{n}{n-\alpha},\infty}}
   \sim
   \|T_{\Omega}^{\alpha}\|_{L^1\rightarrow L^{\frac{n}{n-\alpha},\infty}}
   \sim
   \|T_{\Omega}^{\alpha}\|_{L^1(l^r)\rightarrow L^{\frac{n}{n-\alpha},\infty}(l^r)}
  $;
  \item
   $\displaystyle\lim_{t\rightarrow 0^{+}}
   \Big\|T^{\alpha}_{\Omega}f_t(\cdot)-\frac{\Omega(\cdot)}{|\cdot|^{n-\alpha}}\|f\|_{L^1(\rn)}\Big\|_{L^{\frac{n}{n-\alpha},\infty}(\mathbb{R}^n\backslash B(0,\rho))}=0$
    for every $\rho>0$;
  \item
    $\displaystyle\lim_{t\rightarrow 0^{+}}
    \bigg\|
     \Big\|
      \big\{
       T^{\alpha}_{\Omega}f_{j,t}(\cdot)-\frac{|\Omega(\cdot)|}{|\cdot|^{n-\alpha}}\|f_j\|_{L^1(\rn)}\big
      \}_{j\in\bbN}
     \Big\|_{l^r}
    \bigg\|_{L^{\frac{n}{n-\alpha},\infty}(\mathbb{R}^n\backslash B(0,\rho))}
    =0$
    for every $\rho>0$.
\end{enumerate}
\end{theorem}

As mentioned above, in the scalar-valued setting, our main theorems improve the limiting weak-type behaviors of fractional maximal
and factional integrals in the previous works.
We also remark that, in the vector-valued setting, our main theorems essentially improve \cite[Thm 1.1, Thm 1.2]{HGW-JMAA} with $\al\in (0,n)$.
In fact, without regularity assumption on $\Om$, we give stronger conclusions.
Let us make some explanation for $M_{\Omega}^{\alpha}$, the same argument works for $T_{\Omega}^{\alpha}$.
By Proposition \ref{proposition, comparison for limiting behavior}, if the conclusion (d) in Theorem \ref{thm, maximal} is valid,
we deduce that
\be
     \lim_{t\rightarrow 0^{+}}
     \left|\left\{x\in \rn:
      \Big\|\big\{
       M^{\alpha}_{\Omega}f_{j,t}(\cdot)-\frac{|\Omega(\cdot)|}{|\cdot|^{n-\alpha}}\|f_j\|_{L^1(\rn)}\big
      \}_{j\in\bbN}
     \Big\|_{l^r}>\lambda
    \right\}\right|=0.
\ee
From this and the triangle inequality of $l^r$, we recapture the main conclusion in \cite[Thm 1.1, Thm 1.2]{HGW-JMAA} with $\al\in (0,n)$:
\be
\begin{split}
  \lim_{t\rightarrow 0^{+}}
     \left|\left\{x\in \rn:
      \left|\Big\|\big\{
       M^{\alpha}_{\Omega}f_{j,t}(\cdot)\big\}\Big\|_{l^r}-\Big\|\frac{|\Omega(\cdot)|}{|\cdot|^{n-\alpha}}\big\{\|f_j\|_{L^1(\rn)}\big
      \}_{j\in\bbN}
     \Big\|_{l^r}\right|>\lambda
    \right\}\right|=0.
\end{split}
\ee

The proofs of Theorem \ref{thm, maximal} and \ref{thm, fractional} will be presented in Section 2.

In this paper,
by $C$ we denote  a {positive constant} which
is independent of the main parameters, but it may vary from line to
line. The symbol $f\lesssim g$ represents that $f\leq Cg$ for some
positive constant $C$. If $f\lesssim g$ and $g\lesssim f$,
we then write $f\sim g$.

\section{Proofs of Theorem \ref{thm, maximal} and \ref{thm, fractional}.}
In this section, we present the proofs of our main theorems. We point out that
reduction method plays an important role in our proofs, not only for the simplification of proofs but also
for the improvements of conclusions.

{\it Proof of Theorem \ref{thm, maximal}.}\quad
This proof is divided into several steps.

\textbf{Step 1: $(1)\Longrightarrow (3)\Longrightarrow (5)\Longrightarrow (4)\Longrightarrow (2)$.}
If $\Omega\in L^{\frac{n}{n-\alpha}}(\mathbb{S}^{n-1})$, by a direct calculation we have
\be
\|\Omega\|^{\frac{n}{n-\alpha}}_{L^{\frac{n}{n-\alpha}}(\mathbb{S}^{n-1})}
  =n\lambda^{\frac{n}{n-\alpha}}\Big|\Big\{x\in \mathbb{R}^n: \frac{|\Omega(x)|}{|x|^{n-\alpha}}>\lambda\Big\}\Big|
  =n\Big|\Big\{x\in \mathbb{R}^n: \frac{|\Omega(x)|}{|x|^{n-\alpha}}>1\Big\}\Big|.
\ee
Then we have
$\Big\|\frac{\Omega(\cdot)}{|\cdot|^{n-\alpha}}\Big\|_{L^{\frac{n}{n-\alpha},\infty}}
\sim \|\Omega\|_{L^{\frac{n}{n-\alpha}}(\mathbb{S}^{n-1})}$.
This and the weak type of Young's inequality (see \cite[Theorem 1.2.13]{Grafakos_Classical_2008}) yield that
\be
\|T^{\alpha}_{|\Om|}f\|_{L^{\frac{n}{n-\alpha},\infty}}
=
\left\|\frac{|\Omega(\cdot)|}{|\cdot|^{n-\alpha}}\ast f\right\|_{L^{\frac{n}{n-\alpha},\infty}}
\lesssim
\left\|\frac{\Omega(\cdot)}{|\cdot|^{n-\alpha}}\right\|_{L^{\frac{n}{n-\alpha},\infty}}\|f\|_{L^1(\rn)}
\sim
\|\Omega\|_{L^{\frac{n}{n-\alpha}}(\mathbb{S}^{n-1})}\|f\|_{L^1(\rn)}.
\ee

Noticing that  $T^{\alpha}_{|\Om|}$ is a positive linear operator, if $T^{\alpha}_{|\Om|}$ is bounded from $L^1$ to $L^{\frac{n}{n-\alpha},\infty}$,
then by the vector-valued inequality (see \cite[Proposition 5.5.10]{Grafakos_Classical_2008}),
we have $T_{|\Omega|}^{\alpha}$ is bounded from $L^1(l^r)$ to $L^{\frac{n}{n-\alpha},\infty}(l^r)$
and
$\|T^{\alpha}_{|\Om|}\|_{L^1(l^r)\rightarrow L^{\frac{n}{n-\alpha},\infty}(l^r)}
=\|T^{\alpha}_{|\Om|}\|_{L^1\rightarrow L^{\frac{n}{n-\alpha},\infty}}$.

While if $T_{|\Omega|}^{\alpha}$ is bounded from $L^1(l^r)$ to $L^{\frac{n}{n-\alpha},\infty}(l^r)$,
then the vector-valued boundedness of $M_{\Om}^{\al}$ follows by
\be
\bigg\|
 \Big(
  \sum_{j\in\bbN}|M^{\alpha}_{\Om}f_j|^r
 \Big)^{1/r}
\bigg\|_{L^{\frac{n}{n-\alpha},\infty}}
\leq
\bigg\|
 \Big(
  \sum_{j\in\bbN}|T^{\alpha}_{|\Om|}|f_j||^r
 \Big)^{1/r}
\bigg\|_{L^{\frac{n}{n-\alpha},\infty}}
\lesssim
\bigg\|
 \Big(
  \sum_{j\in\bbN}|f_j|^r
 \Big)^{1/r}
\bigg\|_{L^1}.
\ee

Let $f_j=0$ for $j\geq2$, we obtain the boundedness of $M_{\Om}^{\al}:L^1\rightarrow L^{\frac{n}{n-\alpha},\infty}$
and
$\|M^{\alpha}_{\Om}\|_{L^1\rightarrow L^{\frac{n}{n-\alpha},\infty}}
\leq
 \|M^{\alpha}_{\Om}\|_{L^1(l^r)\rightarrow L^{\frac{n}{n-\alpha},\infty}(l^r)}$.

\textbf{Step 2: $(2)\Longrightarrow (1)$.}
In this part, we will prove that
\ben\label{pf of thm1, 1}
\|\Omega\|_{L^{\frac{n}{n-\alpha}}(\mathbb{S}^{n-1})}
\lesssim \|M_{\Omega}^{\alpha}\|_{L^1\rightarrow L^{\frac{n}{n-\alpha},\infty}}.
\een
In order to achieve this goal, some reduction arguments will be applied first as following.

\textbf{2.1: First reduction.}
In order to verify \eqref{pf of thm1, 1}, we only need to consider the case $\Om\in L^{\fy}(\mathbb{S}^{n-1})$.
In fact, denote $\Om_N(x)=\Om(x)\chi_{\{x:|\Om(x)|\leq N\}}$. If \eqref{pf of thm1, 1} holds for all $\Om_N$, we have
\be
\|\Om\|_{L^{\frac{n}{n-\alpha}}(\mathbb{S}^{n-1})}
=\lim_{N\rightarrow \fy}\|\Om_N\|_{L^{\frac{n}{n-\alpha}}(\mathbb{S}^{n-1})}
\lesssim
\lim_{N\rightarrow \fy}\|M_{\Omega_N}^{\alpha}\|_{L^1\rightarrow L^{\frac{n}{n-\alpha},\infty}}
\leq
\|M_{\Omega}^{\alpha}\|_{L^1\rightarrow L^{\frac{n}{n-\alpha},\infty}}.
\ee

\textbf{2.2: Second reduction.}
In order to verify \eqref{pf of thm1, 1} for $\Om\in L^{\fy}(\mathbb{S}^{n-1})$,
we only need to consider the case $\Om\in Lip(\mathbb{S}^{n-1})$.
Observe that $Lip(\mathbb{S}^{n-1})$ is dense in $L^{\frac{n}{n-\alpha}}(\mathbb{S}^{n-1})$.
For $\Om\in L^{\fy}(\mathbb{S}^{n-1})\subset L^{\frac{n}{n-\alpha}}(\mathbb{S}^{n-1})$,
there exists a sequence of $Lip(\mathbb{S}^{n-1})$ functions,
denoted by $\{\Om_{\ep}\}_{\ep>0}$ such that $\lim_{\ep\rightarrow 0^+}\|\Om_{\ep}-\Om\|_{L^{\frac{n}{n-\alpha}}(\mathbb{S}^{n-1})}=0$.
And we have $\Om-\Om_{\ep}\in L^{\infty}\subset L^{\frac{n}{n-\alpha}}(\mathbb{S}^{n-1})$. 
If \eqref{pf of thm1, 1} holds for all $\Om_{\ep}$, by the quasi triangle inequality we have
\be
\begin{split}
  \|\Om\|_{L^{\frac{n}{n-\alpha}}(\mathbb{S}^{n-1})}
\lesssim &
\|\Om_{\ep}\|_{L^{\frac{n}{n-\alpha}}(\mathbb{S}^{n-1})}+\|\Om-\Om_{\ep}\|_{L^{\frac{n}{n-\alpha}}(\mathbb{S}^{n-1})}
\\
\lesssim &
\|M_{\Omega_{\ep}}^{\alpha}\|_{L^1\rightarrow L^{\frac{n}{n-\alpha},\infty}}+\|\Om-\Om_{\ep}\|_{L^{\frac{n}{n-\alpha}}(\mathbb{S}^{n-1})}
\\
\lesssim &
\|M_{\Omega}^{\alpha}+M_{\Om-\Omega_{\ep}}^{\alpha}\|_{L^1\rightarrow L^{\frac{n}{n-\alpha},\infty}}+\|\Om-\Om_{\ep}\|_{L^{\frac{n}{n-\alpha}}(\mathbb{S}^{n-1})}
\\
\lesssim &
\|M_{\Omega}^{\alpha}\|_{L^1\rightarrow L^{\frac{n}{n-\alpha},\infty}}+\|M_{\Om-\Omega_{\ep}}^{\alpha}\|_{L^1\rightarrow L^{\frac{n}{n-\alpha},\infty}}+\|\Om-\Om_{\ep}\|_{L^{\frac{n}{n-\alpha}}(\mathbb{S}^{n-1})}
\\
\lesssim &
\|M_{\Omega}^{\alpha}\|_{L^1\rightarrow L^{\frac{n}{n-\alpha},\infty}}+\|\Om-\Om_{\ep}\|_{L^{\frac{n}{n-\alpha}}(\mathbb{S}^{n-1})}
+\|\Om-\Om_{\ep}\|_{L^{\frac{n}{n-\alpha}}(\mathbb{S}^{n-1})},
\end{split}
\ee
where in the last inequality we use
$\|M_{\Om-\Omega_{\ep}}^{\alpha}\|_{L^1\rightarrow L^{\frac{n}{n-\alpha},\infty}}\lesssim \|\Om-\Om_{\ep}\|_{L^{\frac{n}{n-\alpha}}(\mathbb{S}^{n-1})}$ established in Step 1.
Letting $\ep\rightarrow 0^+$, we conclude that
\be
\|\Om\|_{L^{\frac{n}{n-\alpha}}(\mathbb{S}^{n-1})}\lesssim \|M_{\Omega}^{\alpha}\|_{L^1\rightarrow L^{\frac{n}{n-\alpha},\infty}}.
\ee

\textbf{2.3: Weak limit for $\Om\in Lip(\mathbb{S}^{n-1})$.} Take $f$ to be a $C_c^{\fy}(\rn)$ function supported on $B(0,1)$,
satisfying $\|f\|_{L^1(\rn)}=1$. We claim that for every $\r>0$,
\ben\label{pf of thm1, 2}
\lim_{t\rightarrow 0^{+}} \Big\|M^{\alpha}_{\Omega}f_t(\cdot)-\frac{|\Omega(\cdot)|}{|\cdot|^{n-\alpha}}\Big\|_{L^{\frac{n}{n-\alpha},\infty}(\mathbb{R}^n\backslash B(0,\rho))}=0.
\een
For $x\in \rn\bs B(0,\r)$ and sufficiently small $t$, we have
\be
  \begin{split}
    M^{\alpha}_{\Omega}f_t(x)
        = &
    \sup_{r>0}\frac{1}{r^{n-\alpha}}\int_{B(x,r)\cap B(0,t)}|\Omega(x-y)||f_t(y)|dy
    \\
    = &
    \sup_{|x|-t\leq r\leq |x|+t}\frac{1}{r^{n-\alpha}}\int_{B(x,r)\cap B(0,t)}|\Omega(x-y)||f_t(y)|dy
    \\
    \in &
    \left[\frac{1}{(|x|+t)^{n-\alpha}}\int_{B(0,t)}|\Omega(x-y)||f_t(y)|dy,\  \frac{1}{(|x|-t)^{n-\alpha}}\int_{B(0,t)}|\Omega(x-y)||f_t(y)|dy\right]
    \\\
    \subset &
    \left[\frac{1-\b_t}{|x|^{n-\alpha}}\int_{B(0,t)}|\Omega(x-y)||f_t(y)|dy,\  \frac{1+\b_t}{|x|^{n-\alpha}}\int_{B(0,t)}|\Omega(x-y)||f_t(y)|dy\right],
  \end{split}
\ee
where $\b_t\rightarrow 0^+$ as $t\rightarrow 0^+$, actually, we could take $\b_t=\r^{n-\a}(\f{1}{(\r-t)^{n-\a}}-\f{1}{(\r+t)^{n-\a}})$.
Write
\be
\frac{|\Om(x)|}{|x|^{n-\al}}=\frac{1}{|x|^{n-\alpha}}\int_{B(0,t)}|\Omega(x)||f_t(y)|dy.
\ee
For $x\in \rn\bs B(0,\r)$ we have
\be
\begin{split}
  M^{\alpha}_{\Omega}f_t(x)-\frac{|\Om(x)|}{|x|^{n-\al}}
  \geq
  \frac{1-\b_t}{|x|^{n-\alpha}}\int_{B(0,t)}(|\Omega(x-y)|-|\Omega(x)|)|f_t(y)|dy-\b_t\frac{|\Om(x)|}{|x|^{n-\al}},
\end{split}
\ee
and
\be
\begin{split}
  M^{\alpha}_{\Omega}f_t(x)-\frac{|\Om(x)|}{|x|^{n-\al}}
  \leq
  \frac{1+\b_t}{|x|^{n-\alpha}}\int_{B(0,t)}(|\Omega(x-y)|-|\Omega(x)|)|f_t(y)|dy+\b_t\frac{|\Om(x)|}{|x|^{n-\al}}.
\end{split}
\ee
From the above two estimates, for $x\in \rn\bs B(0,\r)$ we conclude that
\be
\begin{split}
|M^{\alpha}_{\Omega}f_t(x)-\frac{|\Om(x)|}{|x|^{n-\al}}|
\leq &
\frac{1+\b_t}{|x|^{n-\alpha}}\int_{B(0,t)}(|\Omega(x-y)-\Omega(x)|)|f_t(y)|dy+\b_t\frac{|\Om(x)|}{|x|^{n-\al}}.
\end{split}
\ee
Recalling that $\Om\in Lip(\mathbb{S}^{n-1})$ in this case, for $x\in \rn\bs B(0,\r)$ and $y\in B(0,t)$ with sufficiently small $t$, we obtain
\be
|\Om(x-y)-\Om(x)|
=
|\Om(\frac{x-y}{|x-y|})-\Om(\frac{x}{|x|})|
\lesssim \left|\frac{x-y}{|x-y|}-\frac{x}{|x|}\right|\lesssim \frac{|y|}{|x|}.
\ee
Using the above estimate and the fact $\|f_t\|_{L^1(\rn)}=1$, we continue the estimate by
\be
\begin{split}
|M^{\alpha}_{\Omega}f_t(x)-\frac{|\Om(x)|}{|x|^{n-\al}}|
\lesssim
\frac{1+\b_t}{|x|^{n-\alpha+1}}\int_{B(0,t)}|y||f_t(y)|dy+\b_t\frac{|\Om(x)|}{|x|^{n-\al}}
\lesssim
\frac{(1+\b_t)t}{|x|^{n-\alpha}}+\b_t\frac{|\Om(x)|}{|x|^{n-\al}}.
\end{split}
\ee
Hence, the desired claim \eqref{pf of thm1, 2} follows by
\be
\begin{split}
 \Big\|M^{\alpha}_{\Omega}f_t(\cdot)-\frac{|\Omega(\cdot)|}{|\cdot|^{n-\alpha}}\Big\|_{L^{\frac{n}{n-\alpha},\infty}(\mathbb{R}^n\backslash B(0,\rho))}
 \leq &
 \Big\|\frac{C(1+\b_t)t}{|\cdot|^{n-\alpha}}+\b_t\frac{|\Om(\cdot)|}{|\cdot|^{n-\al}}\Big\|_{L^{\frac{n}{n-\alpha},\infty}(\mathbb{R}^n\backslash B(0,\rho))}
 \\
 \lesssim &
 C(1+\b_t)t\Big\|\frac{1}{|\cdot|^{n-\alpha}}\Big\|_{L^{\frac{n}{n-\alpha},\infty}(\rn)}
 +
\b_t\Big\|\frac{\Om(\cdot)}{|\cdot|^{n-\alpha}}\Big\|_{L^{\frac{n}{n-\alpha},\infty}(\rn)}
\\
\lesssim &
(1+\b_t)t+\b_t\rightarrow 0,\ \ \ \ \ (t\rightarrow 0^+),
\end{split}
\ee
where in the last inequality we use the fact $\frac{1}{|x|^{n-\al}}, \frac{|\Om(x)|}{|x|^{n-\al}}\in L^{\frac{n}{n-\alpha},\infty}(\rn)$.

\textbf{2.4: Upper bound for $\Om\in Lip(\mathbb{S}^{n-1})$.}
Write
\be
\begin{split}
  \Big\|\frac{|\Omega(\cdot)|}{|\cdot|^{n-\alpha}}\Big\|_{L^{\frac{n}{n-\alpha},\infty}(\mathbb{R}^n\backslash B(0,\rho))}
  \lesssim &
  \Big\|M^{\alpha}_{\Omega}f_t(\cdot)-\frac{|\Omega(\cdot)|}{|\cdot|^{n-\alpha}}\Big\|_{L^{\frac{n}{n-\alpha},\infty}(\mathbb{R}^n\backslash B(0,\rho))}
  +
  \Big\|M^{\alpha}_{\Omega}f_t(\cdot)\Big\|_{L^{\frac{n}{n-\alpha},\infty}(\mathbb{R}^n\backslash B(0,\rho))}
  \\
  \leq &
  \Big\|M^{\alpha}_{\Omega}f_t(\cdot)-\frac{|\Omega(\cdot)|}{|\cdot|^{n-\alpha}}\Big\|_{L^{\frac{n}{n-\alpha},\infty}(\mathbb{R}^n\backslash B(0,\rho))}
  +
  \Big\|M^{\alpha}_{\Omega}f_t(\cdot)\Big\|_{L^{\frac{n}{n-\alpha},\infty}(\mathbb{R}^n)}
  \\
  \lesssim &
  \Big\|M^{\alpha}_{\Omega}f_t(\cdot)-\frac{|\Omega(\cdot)|}{|\cdot|^{n-\alpha}}\Big\|_{L^{\frac{n}{n-\alpha},\infty}(\mathbb{R}^n\backslash B(0,\rho))}+\|M_{\Omega}^{\alpha}\|_{L^1\rightarrow L^{\frac{n}{n-\alpha},\infty}}.
\end{split}
\ee
Letting $t\rightarrow 0^+$, using the conclusion in Step 2.3 we obtain
\be
\Big\|\frac{|\Omega(\cdot)|}{|\cdot|^{n-\alpha}}\Big\|_{L^{\frac{n}{n-\alpha},\infty}(\mathbb{R}^n\backslash B(0,\rho))}
  \lesssim
  \|M_{\Omega}^{\alpha}\|_{L^1\rightarrow L^{\frac{n}{n-\alpha},\infty}}.
\ee
Finally, letting $\r\rightarrow 0^+$, we conclude that
\be
\begin{split}
\|\Om\|_{L^{\frac{n}{n-\alpha}}(\mathbb{S}^{n-1})}
\sim
  \Big\|\frac{|\Omega(\cdot)|}{|\cdot|^{n-\alpha}}\Big\|_{L^{\frac{n}{n-\alpha},\infty}(\rn)}
  =
  \lim_{\r\rightarrow 0^+}\Big\|\frac{|\Omega(\cdot)|}{|\cdot|^{n-\alpha}}\Big\|_{L^{\frac{n}{n-\alpha},\infty}(\mathbb{R}^n\backslash B(0,\rho))}
  \lesssim
  \|M_{\Omega}^{\alpha}\|_{L^1\rightarrow L^{\frac{n}{n-\alpha},\infty}}.
\end{split}
\ee
This completes the proof for $\Om\in Lip(\mathbb{S}^{n-1})$. The final conclusion for $\Om\in L^1(\mathbb{S}^{n-1})$
follows by the first and second reduction in Step 2.1 and 2.2.

\textbf{Step 3: $(1),(2),(3),(4),(5)\Longrightarrow (a).$}
By the estimates in Step 1 and Step 2, the conclusion (a) follows by
\be
\begin{split}
  \|\Om\|_{L^{\frac{n}{n-\alpha}}(\mathbb{S}^{n-1})}
   &
   \sim
  \Big\|\frac{|\Omega(\cdot)|}{|\cdot|^{n-\alpha}}\Big\|_{L^{\frac{n}{n-\alpha},\infty}(\rn)}
  \lesssim
  \|M^{\alpha}_{\Om}\|_{L^1\rightarrow L^{\frac{n}{n-\alpha},\infty}}
   \leq
  \|M^{\alpha}_{\Om}\|_{L^1(l^r)\rightarrow L^{\frac{n}{n-\alpha},\infty}(l^r)}
   \\
   &
   \leq
   \|T^{\alpha}_{|\Om|}\|_{L^1(l^r)\rightarrow L^{\frac{n}{n-\alpha},\infty}(l^r)}
   =
  \|T^{\alpha}_{|\Om|}\|_{L^1\rightarrow L^{\frac{n}{n-\alpha},\infty}}
   \lesssim
  \|\Om\|_{L^{\frac{n}{n-\alpha}}(\mathbb{S}^{n-1})}.
\end{split}
\ee

\textbf{Step 4: $(1),(2),(3),(4),(5)\Longrightarrow (b).$}

\textbf{4.1: First reduction.}
In order to verify (b) for $\Om\in L^{\frac{n}{n-\alpha}}(\mathbb{S}^{n-1})$,
we only need to consider the case $\Om\in Lip(\mathbb{S}^{n-1})$.
Since $Lip(\mathbb{S}^{n-1})$ is dense in $L^{\frac{n}{n-\alpha}}(\mathbb{S}^{n-1})$,
there exists a sequence of $Lip(\mathbb{S}^{n-1})$ functions $\{\Om_{\ep}\}_{\ep>0}$ such that 
$\lim_{\ep\rightarrow 0^+}\|\Om_{\ep}-\Om\|_{L^{\frac{n}{n-\alpha}}(\mathbb{S}^{n-1})}=0$.
Observe that
\be
\begin{split}
&\left|M^{\alpha}_{\Omega}f_t(x)-\|f\|_{L^1(\rn)}\frac{|\Omega(x)|}{|x|^{n-\alpha}}\right|
\\
\leq &
\left|M^{\alpha}_{\Omega_{\ep}}f_t(x)-\|f\|_{L^1(\rn)}\frac{|\Omega_{\ep}(x)|}{|x|^{n-\alpha}}\right|
+
\left|M^{\alpha}_{\Omega-\Om_{\ep}}f_t(x)\right|+\|f\|_{L^1(\rn)}\left|\frac{|\Omega(x)-\Om_{\ep}(x)|}{|x|^{n-\alpha}}\right|.
\end{split}
\ee
If (b) holds for all $\Om_{\ep}$, we have
\be
\begin{split}
  &\Big\|M^{\alpha}_{\Omega}f_t(\cdot)-\|f\|_{L^1(\rn)}\frac{|\Omega(\cdot)|}{|\cdot|^{n-\alpha}}\Big\|_{L^{\frac{n}{n-\alpha},\infty}(\mathbb{R}^n\backslash B(0,\rho))}
  \\
  \lesssim &
  \Big\|M^{\alpha}_{\Omega_{\ep}}f_t(\cdot)-\|f\|_{L^1(\rn)}\frac{|\Omega_{\ep}(\cdot)|}{|\cdot|^{n-\alpha}}\Big\|_{L^{\frac{n}{n-\alpha},\infty}(\mathbb{R}^n\backslash B(0,\rho))}
  +
  \Big\|M^{\alpha}_{\Omega-\Om_{\ep}}f_t(\cdot)\Big\|_{L^{\frac{n}{n-\alpha},\infty}(\mathbb{R}^n\backslash B(0,\rho))}
  \\
  + &
  \|f\|_{L^1(\rn)}\Big\|\frac{|\Omega(\cdot)-\Om_{\ep}(\cdot)|}{|\cdot|^{n-\alpha}}\Big\|_{L^{\frac{n}{n-\alpha},\infty}(\mathbb{R}^n\backslash B(0,\rho))}
  \\
  \lesssim &
  \Big\|M^{\alpha}_{\Omega_{\ep}}f_t(\cdot)-\|f\|_{L^1(\rn)}\frac{|\Omega_{\ep}(\cdot)|}{|\cdot|^{n-\alpha}}\Big\|_{L^{\frac{n}{n-\alpha},\infty}(\mathbb{R}^n\backslash B(0,\rho))}
  +
  \|f\|_{L^1(\rn)}\|\Om-\Om_{\ep}\|_{L^{\frac{n}{n-\alpha}}(\mathbb{S}^{n-1})}.
\end{split}
\ee
Letting $t\rightarrow 0^+$, we have
\be
\begin{split}
& \varlimsup_{t\rightarrow 0^+}\Big\|M^{\alpha}_{\Omega}f_t(\cdot)-\|f\|_{L^1(\rn)}\frac{|\Omega(\cdot)|}{|\cdot|^{n-\alpha}}\Big\|_{L^{\frac{n}{n-\alpha},\infty}(\mathbb{R}^n\backslash B(0,\rho))}
\\
\lesssim &
\varlimsup_{t\rightarrow 0^+}
\Big\|M^{\alpha}_{\Omega_{\ep}}f_t(\cdot)-\|f\|_{L^1(\rn)}\frac{|\Omega_{\ep}(\cdot)|}{|\cdot|^{n-\alpha}}\Big\|_{L^{\frac{n}{n-\alpha},\infty}(\mathbb{R}^n\backslash B(0,\rho))}
+\|f\|_{L^1(\rn)}\|\Om-\Om_{\ep}\|_{L^{\frac{n}{n-\alpha}}(\mathbb{S}^{n-1})}
\\
= &
\|f\|_{L^1(\rn)}\|\Om-\Om_{\ep}\|_{L^{\frac{n}{n-\alpha}}(\mathbb{S}^{n-1})}.
\end{split}
\ee
Then, the conclusion follows by letting $\ep\rightarrow 0^+$.

\textbf{4.2: Second reduction.}
In order to verify (b) for $\Om\in Lip(\mathbb{S}^{n-1})$ and $f\in L^1(\rn)$,
we only need to consider the case $f\in C_c^{\fy}(\rn)$.
Recall that $C_c^{\fy}(\rn)$ is dense in $L^1(\rn)$.
We take a sequence of functions $\{f^{(m)}\}_{m\in \mathbb{N}}\subset C_c^{\fy}(\rn)$ such that
$f^{(m)}\rightarrow f$ in $L^1(\rn)$.
Observe that
\be
\begin{split}
&\left|M^{\alpha}_{\Omega}f_t(x)-\|f\|_{L^1(\rn)}\frac{|\Omega(x)|}{|x|^{n-\alpha}}\right|
\\
\leq &
\left|M^{\alpha}_{\Omega}f^{(m)}_t(x)-\|f^{(m)}\|_{L^1(\rn)}\frac{|\Omega(x)|}{|x|^{n-\alpha}}\right|
+
\left|M^{\alpha}_{\Omega}(|f_t-f^{(m)}_t|)(x)\right|
+
\|f-f^{(m)}\|_{L^1(\rn)}\frac{|\Omega(x)|}{|x|^{n-\alpha}}.
\end{split}
\ee
Then, we conclude that
\be
\begin{split}
  &\Big\|M^{\alpha}_{\Omega}f_t(\cdot)-\|f\|_{L^1(\rn)}\frac{|\Omega(\cdot)|}{|\cdot|^{n-\alpha}}\Big\|_{L^{\frac{n}{n-\alpha},\infty}(\mathbb{R}^n\backslash B(0,\rho))}
  \\
  \lesssim &
  \Big\|M^{\alpha}_{\Omega}f^{(m)}_t(\cdot)-\|f^{(m)}\|_{L^1(\rn)}\frac{|\Omega(\cdot)|}{|\cdot|^{n-\alpha}}\Big\|_{L^{\frac{n}{n-\alpha},\infty}(\mathbb{R}^n\backslash B(0,\rho))}
  +
  \Big\|M^{\alpha}_{\Omega}(|f_t-f^{(m)}_t|)(\cdot)\Big\|_{L^{\frac{n}{n-\alpha},\infty}(\mathbb{R}^n\backslash B(0,\rho))}
  \\
  & + 
  \|f-f^{(m)}\|_{L^1(\rn)}\Big\|\frac{|\Omega(\cdot)|}{|\cdot|^{n-\alpha}}\Big\|_{L^{\frac{n}{n-\alpha},\infty}(\mathbb{R}^n\backslash B(0,\rho))}
  \\
  \lesssim &
  \Big\|M^{\alpha}_{\Omega}f^{(m)}_t(\cdot)-\|f^{(m)}\|_{L^1(\rn)}\frac{|\Omega(\cdot)|}{|\cdot|^{n-\alpha}}\Big\|_{L^{\frac{n}{n-\alpha},\infty}(\mathbb{R}^n\backslash B(0,\rho))}
  +\|f-f^{(m)}\|_{L^1(\rn)}.
\end{split}
\ee
If (b) holds for all $f^{(m)}$, letting $t\rightarrow 0^+$ we have
\be
\begin{split}
& \varlimsup_{t\rightarrow 0^+}\Big\|M^{\alpha}_{\Omega}f_t(\cdot)-\|f\|_{L^1(\rn)}\frac{|\Omega(\cdot)|}{|\cdot|^{n-\alpha}}\Big\|_{L^{\frac{n}{n-\alpha},\infty}(\mathbb{R}^n\backslash B(0,\rho))}
\\
\lesssim &
\varlimsup_{t\rightarrow 0^+}
\Big\|M^{\alpha}_{\Omega}f^{(m)}_t(\cdot)-\|f^{(m)}\|_{L^1(\rn)}\frac{|\Omega(\cdot)|}{|\cdot|^{n-\alpha}}\Big\|_{L^{\frac{n}{n-\alpha},\infty}(\mathbb{R}^n\backslash B(0,\rho))}
+
\|f-f^{(m)}\|_{L^1(\rn)}
=
\|f-f^{(m)}\|_{L^1(\rn)}.
\end{split}
\ee
Then, the conclusion follows by letting $m\rightarrow \fy$.

\textbf{4.3: Third reduction.}
In order to verify (b) for $\Om\in Lip(\mathbb{S}^{n-1})$ and $f\in C_c^{\fy}(\rn)$,
we only need to consider the case that $f$ is a $C_c^{\fy}(\rn)$ function supported on $B(0,1)$,
satisfying $\|f\|_{L^1(\rn)}=1$.
In fact, for any nonzero $C_c^{\fy}(\rn)$ function $f$ supported in $B(0,R)$,
then $g(x):=\frac{1}{R^n\|f\|_{L^1(\rn)}}f(x/R)$ is a $C_c^{\fy}(\rn)$ function supported on $B(0,1)$,
satisfying $\|g\|_{L^1(\rn)}=1$.
If (b) holds for $g$, we have
\be
\displaystyle\lim_{t\rightarrow 0^{+}} \Big\|M^{\alpha}_{\Omega}g_t(\cdot)-\frac{|\Omega(\cdot)|}{|\cdot|^{n-\alpha}}\Big\|_{L^{\frac{n}{n-\alpha},\infty}(\mathbb{R}^n\backslash B(0,\rho))}=0.
\ee
Note that $g_t=f_{Rt}/\|f\|_{L^1(\rn)}$. We conclude that
\be
\displaystyle\lim_{t\rightarrow 0^{+}} \Big\|\frac{1}{\|f\|_{L^1(\rn)}}M^{\alpha}_{\Omega}f_t(\cdot)-\frac{|\Omega(\cdot)|}{|\cdot|^{n-\alpha}}\Big\|_{L^{\frac{n}{n-\alpha},\infty}(\mathbb{R}^n\backslash B(0,\rho))}=0,
\ee
which implies that
\be
\displaystyle\lim_{t\rightarrow 0^{+}} \Big\|M^{\alpha}_{\Omega}f_t(\cdot)-\|f\|_{L^1(\rn)}\frac{|\Omega(\cdot)|}{|\cdot|^{n-\alpha}}\Big\|_{L^{\frac{n}{n-\alpha},\infty}(\mathbb{R}^n\backslash B(0,\rho))}=0.
\ee
Using the above reduction arguments, the final conclusion follows by Step 2.3.

\textbf{Step 5: $(1),(2),(3),(4),(5)\Longrightarrow (c).$}
Using the same reduction method as in Step 4, we only need to consider the case $\Om\in Lip(\mathbb{S}^{n-1})$,
and $f$ is a $C_c^{\fy}(\rn)$ function supported on $B(0,1)$, satisfying $\|f\|_{L^1(\rn)}=1$.
In this case, we have
\be
  \begin{split}
    T^{\alpha}_{|\Om|}f_t(x)-\frac{|\Om(x)|}{|x|^{n-\alpha}}
    =
    \int_{B(0,t)}\bigg(\frac{|\Omega(x-y)|}{|x-y|^{n-\al}}-\frac{|\Omega(x)|}{|x|^{n-\al}}\bigg)f_t(y)dy.
  \end{split}
\ee
Since $\Om\in Lip(\mathbb{S}^{n-1})$,
for $x\in \rn\bs B(0,\r)$ and $y\in B(0,t)$ with sufficiently small $t$, we have
\be
\begin{split}
\left|\frac{|\Omega(x-y)|}{|x-y|^{n-\al}}-\frac{|\Omega(x)|}{|x|^{n-\al}}\right|
\lesssim &
\left|\Omega(x-y)\left(\frac{1}{|x-y|^{n-\al}}-\frac{1}{|x|^{n-\al}}\right)\right|
+
\left|\frac{|\Omega(x-y)|-|\Omega(x)|}{|x|^{n-\al}}\right|
\\
\lesssim &
\frac{|y|}{|x|^{n-\al+1}}\lesssim \frac{|y|}{|x|^{n-\al}}.
\end{split}
\ee
Then,
\be
\begin{split}
  \left|T^{\alpha}_{|\Omega|}f_t(x)-\frac{|\Omega(x)|}{|x|^{n-\alpha}}\right|
  \lesssim
  \frac{1}{|x|^{n-\al}}\int_{B(0,t)}|y||f_t(y)|dy\leq\frac{t}{|x|^{n-\al}}.
\end{split}
\ee
The final conclusion follows by
\be
\begin{split}
  \Big\|T^{\alpha}_{|\Omega|}f_t(\cdot)-\frac{|\Omega(\cdot)|}{|\cdot|^{n-\alpha}}\Big\|_{L^{\frac{n}{n-\alpha},\infty}(\mathbb{R}^n\backslash B(0,\rho))}
  \lesssim
  t\||x|^{\al-n}\|_{L^{\frac{n}{n-\alpha},\infty}(\mathbb{R}^n)}\lesssim t,
\end{split}
\ee
which tends to zero as $t\rightarrow 0^+$.

\textbf{Step 6: $(1),(2),(3),(4),(5)\Longrightarrow (d),(e).$}
In fact, the vector-valued case follows directly by the scalar-valued case and a reduction argument.
Denote by $A$ the operator $M^{\alpha}_{\Omega}$ or $T^{\alpha}_{|\Omega|}$.
By the quasi triangle inequality and the boundedness of $A$, we write
\begin{align*}
 &
    \bigg\|
     \Big\|
      \big\{
       Af_{j,t}(\cdot)-\frac{|\Omega(\cdot)|}{|\cdot|^{n-\alpha}}\|f_j\|_{L^1(\rn)}
       \big\}_{j\in\bbN}
     \Big\|_{l^r}
    \bigg\|_{L^{\frac{n}{n-\alpha},\infty}(\mathbb{R}^n\backslash B(0,\rho))}
    \\
    \lesssim
    &
    \bigg\|
     \Big\|
      \big\{
       Af_{j,t}(\cdot)-\frac{|\Omega(\cdot)|}{|\cdot|^{n-\alpha}}\|f_j\|_{L^1(\rn)}
       \big\}_{1\leq j\leq N}
     \Big\|_{l^r}
    \bigg\|_{L^{\frac{n}{n-\alpha},\infty}(\mathbb{R}^n\backslash B(0,\rho))}
    \\
    &+
    \bigg\|
     \Big\|
      \big\{
       Af_{j,t}(\cdot)
       \big\}_{j\geq N}
     \Big\|_{l^r}
    \bigg\|_{L^{\frac{n}{n-\alpha},\infty}(\mathbb{R}^n\backslash B(0,\rho))}
    +
    \bigg\|
     \Big\|
      \big\{
       \frac{|\Omega(\cdot)|}{|\cdot|^{n-\alpha}}\|f_j\|_{L^1(\rn)}
       \big\}_{j\geq N}
     \Big\|_{l^r}
    \bigg\|_{L^{\frac{n}{n-\alpha},\infty}(\mathbb{R}^n\backslash B(0,\rho))}
    \\
    \lesssim
    &
    \sum_{1\leq j\leq N}
     \bigg\|
       Af_{j,t}(\cdot)-\frac{|\Omega(\cdot)|}{|\cdot|^{n-\alpha}}\|f_j\|_{L^1(\rn)}
    \bigg\|_{L^{\frac{n}{n-\alpha},\infty}(\mathbb{R}^n\backslash B(0,\rho))}
    +\bigg\|
      \big\{
       f_{j,t}(\cdot)
       \big\}_{j\geq N}
    \bigg\|_{L^1(l^r)}
    \\
    &+\|\Om\|_{L^{\frac{n}{n-\alpha}}(\mathbb{S}^{n-1})}\cdot
    \Big\|
      \big\{\|f_j\|_{L^1(\rn)}
       \big\}_{j\geq N}
     \Big\|_{l^r}
    \\
    \lesssim &
    \sum_{1\leq j\leq N}
     \bigg\|
       Af_{j,t}(\cdot)-\frac{|\Omega(\cdot)|}{|\cdot|^{n-\alpha}}\|f_j\|_{L^1(\rn)}
    \bigg\|_{L^{\frac{n}{n-\alpha},\infty}(\mathbb{R}^n\backslash B(0,\rho))}
    +
    \bigg\|
      \big\{
       f_{j}(\cdot)
       \big\}_{j\geq N}
    \bigg\|_{L^1(l^r)}=: I_{N,t}+R_N,
\end{align*}
where in the last inequality we use the fact $\Big\|
      \big\{\|f_j\|_{L^1(\rn)}
       \big\}_{j\geq N}
     \Big\|_{l^r}\lesssim \bigg\|
      \big\{
       f_{j}(\cdot)
       \big\}_{j\geq N}
    \bigg\|_{L^1(l^r)}$ by Minknowski's inequality.
Letting $t\rightarrow 0^+$, by the corresponding results of scalar-valued case, we conclude that
\ben\label{pf1}
\varlimsup_{t\rightarrow 0^+}\bigg\|
     \Big\|
      \big\{
       Af_{j,t}(\cdot)-\frac{|\Omega(\cdot)|}{|\cdot|^{n-\alpha}}\|f_j\|_{L^1(\rn)}
       \big\}_{j\in\bbN}
     \Big\|_{l^r}
    \bigg\|_{L^{\frac{n}{n-\alpha},\infty}(\mathbb{R}^n\backslash B(0,\rho))}\lesssim R_N.
\een
Recalling that
$\bigg\|\big\{f_{j}(\cdot)\big\}_{j\in \mathbb{N}}\bigg\|_{L^1(l^r)}<\fy$,
by the Lebesgue dominated convergence theorem we have $\lim_{N\rightarrow \fy}R_N=0$.
Hence, the desired conclusion follows by letting $N\rightarrow \fy$ in \eqref{pf1}.

\medskip

{\it \textbf{Proof of Theorem \ref{thm, fractional}.}}\quad
By Theorem \ref{thm, maximal} and $\Om\in L^{\frac{n}{n-\alpha}}(\mathbb{S}^{n-1})$,
$T_{|\Omega|}^{\alpha}$ is bounded from $L^1$ to $L^{\frac{n}{n-\alpha},\infty}$ with
$\|\Omega \|_{L^{\frac{n}{n-\alpha}}(\mathbb{S}^{n-1})}
\sim
\|T_{|\Omega|}^{\alpha}\|_{L^1\rightarrow L^{\frac{n}{n-\alpha},\infty}}
\sim
\|T_{|\Omega|}^{\alpha}\|_{L^1(l^r)\rightarrow L^{\frac{n}{n-\alpha},\infty}(l^r)}
$.
The boundedness of $T_{\Om}^{\al}$ follows by
\be
\|T_{\Omega}^{\alpha}\|_{L^1\rightarrow L^{\frac{n}{n-\alpha},\infty}}
\leq
\|T_{\Omega}^{\alpha}\|_{L^1(l^r)\rightarrow L^{\frac{n}{n-\alpha},\infty}(l^r)}
\leq
\|T_{|\Omega|}^{\alpha}\|_{L^1(l^r)\rightarrow L^{\frac{n}{n-\alpha},\infty}(l^r)}
\sim
\|\Omega \|_{L^{\frac{n}{n-\alpha}}(\mathbb{S}^{n-1})}.
\ee

\textbf{Proof for (1):}
We only need to prove that
\ben\label{pf of thm1, 3}
\|\Omega\|_{L^{\frac{n}{n-\alpha}}(\mathbb{S}^{n-1})}
\lesssim \|T_{\Omega}^{\alpha}\|_{L^1\rightarrow L^{\frac{n}{n-\alpha},\infty}}.
\een
By a similar reduction argument as in Step 2.2 in the proof of Theorem \ref{thm, maximal},
we only need to consider the case $\Om\in Lip(\mathbb{S}^{n-1})$.
In fact, we can take a sequence of $Lip(\mathbb{S}^{n-1})$ functions,
denoted by $\{\Om_{\ep}\}_{\ep>0}$ such that $\lim_{\ep\rightarrow 0^+}\|\Om_{\ep}-\Om\|_{L^{\frac{n}{n-\alpha}}(\mathbb{S}^{n-1})}=0$.
If \eqref{pf of thm1, 3} holds for all $\Om_{\ep}$, we have
\be
\begin{split}
  \|\Om\|_{L^{\frac{n}{n-\alpha}}(\mathbb{S}^{n-1})}
\leq &
\|\Om_{\ep}\|_{L^{\frac{n}{n-\alpha}}(\mathbb{S}^{n-1})}+\|\Om-\Om_{\ep}\|_{L^{\frac{n}{n-\alpha}}(\mathbb{S}^{n-1})}
\\
\leq &
C\|T_{\Omega_{\ep}}^{\alpha}\|_{L^1\rightarrow L^{\frac{n}{n-\alpha},\infty}}+\|\Om-\Om_{\ep}\|_{L^{\frac{n}{n-\alpha}}(\mathbb{S}^{n-1})}
\\
= &
C\|T_{\Omega}^{\alpha}+T_{\Om-\Omega_{\ep}}^{\alpha}\|_{L^1\rightarrow L^{\frac{n}{n-\alpha},\infty}}+\|\Om-\Om_{\ep}\|_{L^{\frac{n}{n-\alpha}}(\mathbb{S}^{n-1})}
\\
\leq &
C\|T_{\Omega}^{\alpha}\|_{L^1\rightarrow L^{\frac{n}{n-\alpha},\infty}}+C\|T_{\Om-\Omega_{\ep}}^{\alpha}\|_{L^1\rightarrow L^{\frac{n}{n-\alpha},\infty}}+\|\Om-\Om_{\ep}\|_{L^{\frac{n}{n-\alpha}}(\mathbb{S}^{n-1})}
\\
\leq &
C\|T_{\Omega}^{\alpha}\|_{L^1\rightarrow L^{\frac{n}{n-\alpha},\infty}}+C\|\Om-\Om_{\ep}\|_{L^{\frac{n}{n-\alpha}}(\mathbb{S}^{n-1})}
+\|\Om-\Om_{\ep}\|_{L^{\frac{n}{n-\alpha}}(\mathbb{S}^{n-1})}.
\end{split}
\ee
Letting $\ep\rightarrow 0^+$, we conclude that
\be
\|\Om\|_{L^{\frac{n}{n-\alpha}}(\mathbb{S}^{n-1})}\lesssim \|T_{\Omega}^{\alpha}\|_{L^1\rightarrow L^{\frac{n}{n-\alpha},\infty}}.
\ee

For $\Om\in Lip(\mathbb{S}^{n-1})$, we take $f$ to be a $C_c^{\fy}(\rn)$ function supported on $B(0,1)$,
satisfying $\|f\|_{L^1(\rn)}=1$.For every $\r>0$, we will verify that
\ben\label{pf of thm1, 4}
\lim_{t\rightarrow 0^{+}} \Big\|T^{\alpha}_{\Omega}f_t(\cdot)-\frac{\Omega(\cdot)}{|\cdot|^{n-\alpha}}\Big\|_{L^{\frac{n}{n-\alpha},\infty}(\mathbb{R}^n\backslash B(0,\rho))}=0.
\een
Write
\be
  \begin{split}
    T^{\alpha}_{\Om}f_t(x)-\frac{\Om(x)}{|x|^{n-\alpha}}
    =
    \int_{B(0,t)}\bigg(\frac{\Omega(x-y)}{|x-y|^{n-\al}}-\frac{\Omega(x)}{|x|^{n-\al}}\bigg)f_t(y)dy.
  \end{split}
\ee
For $x\in \rn\bs B(0,\r)$ and $y\in B(0,t)$ for sufficiently small $t$, we have
\be
\begin{split}
\left|\frac{\Omega(x-y)}{|x-y|^{n-\al}}-\frac{\Omega(x)}{|x|^{n-\al}}\right|
\lesssim
\frac{|y|}{|x|^{n-\al+1}}\lesssim \frac{|y|}{|x|^{n-\al}}.
\end{split}
\ee
Then,
\be
\begin{split}
  \left|T^{\alpha}_{\Omega}f_t(x)-\frac{\Omega(x)}{|x|^{n-\alpha}}\right|
  \lesssim
  \frac{1}{|x|^{n-\al}}\int_{B(0,t)}|y||f_t(y)|dy\leq\frac{t}{|x|^{n-\al}}.
\end{split}
\ee
The final conclusion follows by
\be
\begin{split}
  \Big\|T^{\alpha}_{\Omega}f_t(\cdot)-\frac{\Omega(\cdot)}{|\cdot|^{n-\alpha}}\Big\|_{L^{\frac{n}{n-\alpha},\infty}(\mathbb{R}^n\backslash B(0,\rho))}
  \lesssim
  t\||x|^{\al-n}\|_{L^{\frac{n}{n-\alpha},\infty}(\mathbb{R}^n)}\lesssim t,
\end{split}
\ee
which tends to zero as $t\rightarrow 0^+$.

From this, letting $t\rightarrow 0^+$ in the following estimate:
\be
\begin{split}
  \Big\|\frac{\Omega(\cdot)}{|\cdot|^{n-\alpha}}\Big\|_{L^{\frac{n}{n-\alpha},\infty}(\mathbb{R}^n\backslash B(0,\rho))}
  \lesssim &
  \Big\|T^{\alpha}_{\Omega}f_t(\cdot)-\frac{\Omega(\cdot)}{|\cdot|^{n-\alpha}}\Big\|_{L^{\frac{n}{n-\alpha},\infty}(\mathbb{R}^n\backslash B(0,\rho))}
  +
  \Big\|T^{\alpha}_{\Omega}f_t(\cdot)\Big\|_{L^{\frac{n}{n-\alpha},\infty}(\mathbb{R}^n\backslash B(0,\rho))}
  \\
  \lesssim &
  \Big\|T^{\alpha}_{\Omega}f_t(\cdot)-\frac{\Omega(\cdot)}{|\cdot|^{n-\alpha}}\Big\|_{L^{\frac{n}{n-\alpha},\infty}(\mathbb{R}^n\backslash B(0,\rho))}+\|T_{\Omega}^{\alpha}\|_{L^1\rightarrow L^{\frac{n}{n-\alpha},\infty}},
\end{split}
\ee
we obtain
\be
\Big\|\frac{\Omega(\cdot)}{|\cdot|^{n-\alpha}}\Big\|_{L^{\frac{n}{n-\alpha},\infty}(\mathbb{R}^n\backslash B(0,\rho))}
  \lesssim
  \|T_{\Omega}^{\alpha}\|_{L^1\rightarrow L^{\frac{n}{n-\alpha},\infty}}.
\ee
Finally, by letting $\r\rightarrow 0^+$, we conclude that
\be
\begin{split}
\|\Om\|_{L^{\frac{n}{n-\alpha}}(\mathbb{S}^{n-1})}
\sim
  \Big\|\frac{\Omega(\cdot)}{|\cdot|^{n-\alpha}}\Big\|_{L^{\frac{n}{n-\alpha},\infty}(\rn)}
  =
  \lim_{\r\rightarrow 0^+}\Big\|\frac{\Omega(\cdot)}{|\cdot|^{n-\alpha}}\Big\|_{L^{\frac{n}{n-\alpha},\infty}(\mathbb{R}^n\backslash B(0,\rho))}
  \lesssim
  \|T_{\Omega}^{\alpha}\|_{L^1\rightarrow L^{\frac{n}{n-\alpha},\infty}}.
\end{split}
\ee

\textbf{Proof for (2)(3):}
Using the same reduction method as in Step 4 in the proof of Theorem \ref{thm, maximal}, we only need to consider the case $\Om\in Lip(\mathbb{S}^{n-1})$,
and $f$ is a $C_c^{\fy}(\rn)$ function supported on $B(0,1)$, satisfying $\|f\|_{L^1(\rn)}=1$.
In this case, we have
\be
  \begin{split}
    T^{\alpha}_{\Om}f_t(x)-\frac{\Om(x)}{|x|^{n-\alpha}}
    =
    \int_{B(0,t)}\bigg(\frac{\Omega(x-y)}{|x-y|^{n-\al}}-\frac{\Omega(x)}{|x|^{n-\al}}\bigg)f_t(y)dy.
  \end{split}
\ee
As in Step 5 in the proof of Theorem $\ref{thm, maximal}$,
for $x\in \rn\bs B(0,\r)$ and $y\in B(0,t)$ for sufficiently small $t$, we have
\be
\begin{split}
\left|\frac{\Omega(x-y)}{|x-y|^{n-\al}}-\frac{\Omega(x)}{|x|^{n-\al}}\right|
\lesssim \frac{|y|}{|x|^{n-\al}}.
\end{split}
\ee
Then,
\be
\begin{split}
  \left|T^{\alpha}_{\Omega}f_t(x)-\frac{\Omega(x)}{|x|^{n-\alpha}}\right|
  \lesssim
  \frac{1}{|x|^{n-\al}}\int_{B(0,t)}|y||f_t(y)|dy\leq\frac{t}{|x|^{n-\al}}.
\end{split}
\ee
The final conclusion follows by
\be
\begin{split}
  \Big\|T^{\alpha}_{\Omega}f_t(\cdot)-\frac{\Omega(\cdot)}{|\cdot|^{n-\alpha}}\Big\|_{L^{\frac{n}{n-\alpha},\infty}(\mathbb{R}^n\backslash B(0,\rho))}
  \lesssim
  t\||x|^{\al-n}\|_{L^{\frac{n}{n-\alpha},\infty}(\mathbb{R}^n)}\lesssim t,
\end{split}
\ee
which tends to zero as $t\rightarrow 0^+$.
The vector-valued case $(3)$ follows by the same argument as in Step 6 in the proof of Theorem \ref{thm, maximal}.

\begin{remark}
  In order to drop the smoothness assumption of $\Om$,
  the key point is to reduce the conclusion as much as possible before the detailed estimates.
  A similar method may work for other types of operators.
  In general, for an operator with kernel $K$, denoted by $T_K$,
  if $T_K(f)$ is sublinear with respect to $K$ and $f$ respectively,
  and if the following boundedness result holds:
  \be
  \|T_Kf\|_{L^{q,\fy}(\rn)}\leq C\|K\|_{Z(\rn)}\|f\|_{L^1(\rn)}
  \ee
  with $Z\subset L^{q,\fy}$,
  then the set consisting of the function pair $(K,f)$ for which
  the type-1 convergence holds, i.e.,
  \be
  \lim_{t\rightarrow0^+}\|T_K(f_t)-K\|f\|_{L^1(\rn)}\|_{L^{q,\fy}(\rn\bs B(0,\r))}=0,
  \ee
  is closed in the usual topology of $Z(\rn)\oplus L^1(\rn)$.
  Hence, in order to conclude the type-1 convergence, we only need to
  conclude it for that $(K,f)$ belongs to some dense subspace of $Z(\rn)\oplus L^1(\rn)$.
\end{remark}

\end{document}